\documentclass[11pt]{article}
\usepackage[colorlinks=true,
linkcolor=webgreen,
filecolor=webbrown,
citecolor=webgreen]{hyperref}

\usepackage{amssymb,psfig,epsfig}
\usepackage{asa}

\definecolor{webgreen}{rgb}{0,.5,0}
\definecolor{webbrown}{rgb}{.6,0,0}

\newtheorem{lemma}{Lemma}
\newtheorem{theorem}[lemma]{Theorem}

\newcommand{\eqn}[1]{(\ref{#1})}
\newcommand{\eeq}{\end{equation}}
\newcommand{\beql}[1]{\begin{equation}\label{#1}}
\newcommand{\bsq}{{\vrule height .9ex width .8ex depth -.1ex }}
\newcommand{\La}{\Lambda}
\newcommand{\af}{\alpha}
\newcommand{\sV}{{\mathcal V}}
\newcommand{\RR}{{\mathbb R}}

\makeatletter
\def\@sect#1#2#3#4#5#6[#7]#8{\ifnum #2>\c@secnumdepth
     \def\@svsec{}\else
     \refstepcounter{#1}\edef\@svsec{\csname the#1\endcsname.\hskip .75em }\fi
     \@tempskipa #5\relax
      \ifdim \@tempskipa>\z@
        \begingroup #6\relax
          \@hangfrom{\hskip #3\relax\@svsec}{\interlinepenalty \@M #8\par}%
        \endgroup
       \csname #1mark\endcsname{#7}\addcontentsline
         {toc}{#1}{\ifnum #2>\c@secnumdepth \else
                      \protect\numberline{\csname the#1\endcsname}\fi
                    #7}\else
        \def\@svsechd{#6\hskip #3\@svsec #8\csname #1mark\endcsname
                      {#7}\addcontentsline
                           {toc}{#1}{\ifnum #2>\c@secnumdepth \else
                             \protect\numberline{\csname the#1\endcsname}\fi
                       #7}}\fi
     \@xsect{#5}}
\def\@begintheorem#1#2{\it \trivlist \item[\hskip \labelsep{\bf #1\ #2.}]}

\def\section{\@startsection {section}{1}{\z@}{-3.5ex plus -1ex minus 
 -.2ex}{2.3ex plus .2ex}{\normalsize\bf}}
\makeatother
\makeatletter
\def\subsection{\@startsection {subsection}{1}{\z@}{-3.5ex plus -1ex minus
 -.2ex}{2.3ex plus .2ex}{\normalsize\bf}}

\makeatother

\begin{document}
\begin{center}
{\large {\bf The Lattice of $N$-Run Orthogonal Arrays}} \\
\vspace{1.5\baselineskip}
E. M. Rains and N. J. A. Sloane \\
\vspace*{1\baselineskip}
Information Sciences Research Center \\
AT\&T Shannon Lab \\
Florham Park, New Jersey 07932-0971 \\ [+.25in]
John Stufken \\
Department of Statistics \\
Iowa State University \\
Ames, IA 50011 \\
\vspace{1.5\baselineskip}
April 20, 2000 \\
\vspace{1.5\baselineskip}
{\bf ABSTRACT}
\vspace{.5\baselineskip}
\end{center}
\setlength{\baselineskip}{1.5\baselineskip}

If the number of runs in a (mixed-level) orthogonal array of strength 2 is
specified, what numbers of levels and factors are possible?
The collection of possible sets of parameters for orthogonal arrays with $N$
runs has a natural lattice structure,
induced by the ``expansive replacement'' construction method.
In particular the dual atoms in this lattice are the most important
parameter sets, since any other parameter set for an $N$-run orthogonal array can be constructed from them.
To get a sense for the number of dual atoms, and to begin to understand the lattice as a function
of $N$, we investigate the height and the size of the lattice.
It is shown that the height is at most $\lfloor c(N-1) \rfloor$,
where $c= 1.4039 \ldots$,
and that there is an infinite sequence of values of $N$ for which
this bound is attained.
On the other hand, the number of nodes in the lattice is bounded above
by a superpolynomial function of $N$ (and superpolynomial growth does occur for certain
sequences of values of $N$).
Using a new construction based on ``mixed spreads'', all parameter sets with 64 runs are determined.
Four of these 64-run orthogonal arrays appear to be new.

\clearpage
\section{Introduction}
Although mixed-level (or asymmetrical) orthogonal arrays have been 
the subject of a number of papers in recent years (see Chapter 9 of
Hedayat, Sloane and Stufken, 1999,
for references),
it seems fair to say that we know much less about them than about fixed-level orthogonal arrays 
(in which all factors have the same number of levels).
For example, there is no analogue for mixed orthogonal arrays of one of the
most powerful construction methods for fixed-level arrays, that based on linear codes (see Chapters 4 and 5 of Hedayat, Sloane and Stufken, 1999).

Again, there are many instances where the linear programming bound for fixed-level orthogonal arrays gives
the correct answer for the minimal number of runs needed for a specified
number of factors.
There {\em is} a linear programming bound for mixed arrays (Sloane
and Stufken, 1996), but it is less effective than in the fixed-level case --- it
ignores too much of the combinatorial nature of the problem
(especially when the levels involve more than one prime number),
and, though generally stronger than the Rao bound, 
does not give correct answers as often as in the fixed-level case.

A mixed orthogonal array
$OA(N, s_1^{k_1} s_2^{k_2} \cdots s_v^{k_v}, t)$ is an array of size $N \times k$,
where $k= k_1 + k_2 + \cdots + k_v$ is the total number of factors,
in which the first $k_1$ columns have symbols from $\{0,1, \ldots, s_1 -1\}$,
the next $k_2$ columns have symbols from $\{0,1, \ldots, s_2 -1\}$, and so on,
with the property that in any $N \times t$ subarray every possible  $t$-tuple of symbols occurs
an equal number of times as a row.
We usually assume $2 \le s_1 < s_2 < \cdots$ and all $k_i \ge 1$.
Except in Section 5, only arrays of strength 2 will be considered, and
we will usually omit $t$ from
the symbol for the array.

We refer to $(N, s_1^{k_1} s_2^{k_2} \ldots )$ as the
{\em parameter set} for an $OA(N,s_1^{k_1} s_2^{k_2} \ldots )$.
We also allow the parameter set $(N, 1^1 )$, corresponding
to the trivial array consisting of a single column of $N$ 0's.
In this paper we consider the question:
if $N$ is specified, how many different parameter sets are possible?

Given an array $A = OA (N, s_1^{k_1} s_2^{k_2} \ldots )$, other $N$-run
arrays can be obtained from it by the
{\em expansive replacement method}.
Let $S$ be one of the $s_i$ occurring in $A$, and suppose $B$ is an $OA(S,t_1^{l_1} t_2^{l_2} \ldots )$.
The expansive replacement method replaces a single column of $A$ at $S$
levels by the rows of $B$.
For example, if $A= OA( 16, 2^3 4^4 )$ and $B= OA(4, 2^3 )$,
we obtain an $OA(16, 2^6 4^3 )$.
If $B$ is a trivial array $OA(S, 1^1 )$, we are simply deleting one of the $S$-level factors from $A$.
E.g. taking $S=2$,
an $OA( 24, 2^{20} 4^1 )$ trivially produces an $OA(24, 2^{19} 4^1 )$.
The expansive replacement method also includes replacing a factor at $s$ levels by a factor at $s'$ levels, if $s'$ divides $s$.
For further details about the expansive
replacement method see Chapter 9 of Hedayat, Sloane and Stufken, 1999.

Let $A$ and $B$ be parameter sets for orthogonal arrays with $N$ runs.
We say that $B$ is {\em dominated by} $A$
if an orthogonal array with parameter set $B$ can be obtained from an orthogonal array with
parameter set $A$ by a sequence of expansive replacements.

Using ``dominance'' as the relation, the parameter sets for orthogonal arrays with $N$ runs form a partially ordered set, which we denote by $\Lambda_N$
(Hedayat, Sloane and Stufken, 1999 p. 335).

$\Lambda_N$ has a unique maximal element $(N, N^1 )$ (corresponding to the
trivial array with one factor at $N$ levels) and a unique minimal element
$(N, 1^1 )$.
It is straightforward to verify that meet $(\wedge)$ and join $(\vee )$ are
well-defined for this relation (we omit the proof), so $\Lambda_N$ is
in fact a {\em lattice} (cf. Welsh, 1976;
Trotter, 1995).

If an $OA(N, s_1^{k_1} s_2^{k_2} \ldots )$ exists,
then necessarily we must have:
\begin{itemize}
\item[(C1)]
$s_i$ divides $N$, for all $i$,
\item[(C2)]
$s_i^2$ divides $N$, if $k_i \ge 2$,
\item[(C3)]
$s_i s_j$ divides $N$, if $i \neq j$,
\item[(C4)]
the Rao bound holds:
\beql{EqR}
N-1 \ge k_1 (s_1 -1) + k_2 (s_2 -1) + \cdots \,,
\eeq
\item[(C5)]
the linear programming bound holds (see Sloane and Stufken, 1996).
\end{itemize}

These conditions are certainly not sufficient for an array to exist,
and it appears to be difficult to test if an orthogonal array does exist
with a putative parameter set satisfying (C1)--(C5).
A further difficulty is that in order to construct $\Lambda_N$ it is necessary
to know $\Lambda_d$ for all proper divisors $d$ of $N$.

To avoid these difficulties we define a second lattice, the
{\em idealized lattice}
$\Lambda'_N$:
this has as nodes all putative parameter sets satisfying
(C1) to (C4), with the dominance relation as before,
except that in the expansive replacement method we may now make use of any
of the nodes of any $\Lambda'_d$ for $d$ dividing $N$.

Constructing $\Lambda'_N$ is much easier than constructing
$\Lambda_N$, since essentially all we need to do is enumerate the solutions
to \eqn{EqR}.
Of course $\Lambda_N$ is a sublattice of $\Lambda'_N$.

To avoid having to repeat the adjective ``putative'',
from now on we will use ``parameter set'' to mean any symbol $(N, s_1^{k_1} s_2^{k_2} \ldots )$ satisfying conditions (C1) to (C4).
The parameter sets are precisely the nodes of $\Lambda'_N$.
If a parameter set is also a node of $\Lambda_N$ then it is implied that an
$OA(N, s_1^{k_1} s_2^{k_2} \ldots )$ does exist, i.e. that the parameter set
is realized by an orthogonal array.

\begin{figure}[htp]
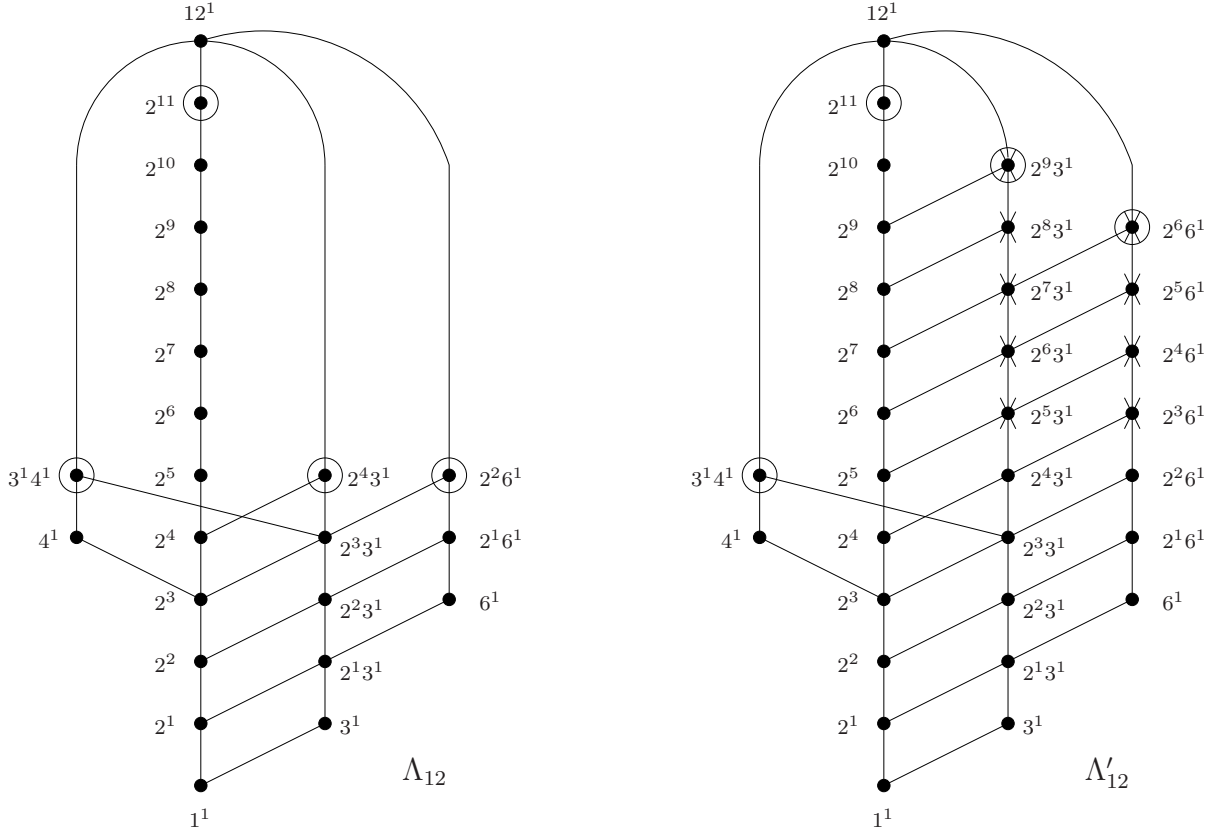

\begin{center}
\input fg1.pstex_t
\end{center}
\caption{(Left) Lattice $\Lambda_{12}$, showing all (true)
parameter sets for 12-run orthogonal arrays.
There are 23 nodes, four dual atoms (circled) and the height is 12.
(Right) Idealized lattice $\Lambda'_{12}$, showing all (putative)
parameter sets satisfying conditions (C1)--(C4).
No arrays exist for the nine nodes marked $\times$.
There are 32 nodes, four dual atoms (circled) and again the height is 12.}
\label{fg1}
\end{figure}

It is convenient to represent $\Lambda_N$ and $\Lambda'_N$ by their Hasse diagrams (cf. Welsh, 1976, p. 45).
These diagrams are drawn ``from the bottom up'',
with $(N, 1^1)$ as the root node at the bottom (Figure \ref{fg1} shows $\La_{12}$ and $\La'_{12}$).
The {\em height} of a parameter set is the number of edges in the longest path
from that node to the root.
A node of height $i$ appears on the $i$th level of the diagram.
The height of the maximal element $(N, N^1)$ will be denoted by $ht(N)$.

The {\em atoms} in $\La_N$ (those nodes just above the root) are
precisely the parameter sets $(N, p^1)$ for the primes $p$ dividing $N$.

The {\em dual atoms} in $\La_N$ (those nodes just below the maximal element) are
especially interesting, since they dominate all other parameter sets.

We can now state our main results.
\begin{theorem}\label{th1}
(i) For all $N$,
\beql{Eq2}
ht(N) \le \lfloor c (N-1) \rfloor ~,
\eeq
where
\beql{Eq2a}
c ~=~ \sum_{i=0}^{\infty} \frac{1}{2^{2^i}-1} ~=~ 1.4039 \ldots ~.
\eeq

(ii) If $N= 2^{2^m}$ $(m \ge 0)$ then $ht(N) = \lfloor c (N-1) \rfloor$.
\end{theorem}

\noindent

Let $T(N)$ (resp. $T' (N)$) denote the total number of nodes in $\La_N$ (resp. $\La'_N$).

\begin{theorem}\label{th2}
If $N=2^{n}$,
\beql{Eq5}
\frac{1}{4} ( \log_2 N)^2 (1+o(1))
~ \le~  \log_2 T(N) ~ \le~  \frac{3}{8}
(\log_2 N)^2 (1+o(1)) \,.
\eeq
\end{theorem}

\begin{theorem}\label{th3}
There is a constant $c_1$ such that for all $N$,
\beql{Eq3}
\ln \ln T(N) ~ \le ~ c_1 ~ \frac{\ln N}{\ln \ln N}  ~ (1+ o(1)) \,.
\eeq
\end{theorem}

\paragraph{Remarks.}
(i) The bounds in \eqn{Eq5} and \eqn{Eq3} also apply to $T' (N)$.

(ii) Theorem \ref{th2} shows that when $N=2^{n}$, $T(N)$ grows very roughly
like $N^{a \log_2 N}$, for some constant $a$ between
$\frac{1}{4}$ and $\frac{3}{8}$.
This is a
``superpolynomial'' function of $N$, meaning that it grows faster than any polynomial in $N$.

(iii) It appears (although we have not proved this) that the upper bound
in \eqn{Eq3} can be achieved by taking $N$ to be a certain product of powers of the first $m$ primes,
where $m$ is about
$$\frac{1}{2e} ~ \frac{\ln N}{\ln \ln N}$$
(see Section 7).
In other words, it appears that there is an infinite sequence of values
of $N$ for which $T(N)$ grows {\em very} roughly like
$$\exp (N^{c_2 / \ln \ln N }) ~,$$
where $c_2$ is a constant.
This is again a superpolynomial function of $N$, and is now
close to being an exponential function, since $\ln \ln N$ grows slowly.

The above discussion has shown that there is an infinite sequence of values of $N$ for which 
the number of nodes in $\La_N$ grows superpolynomially, while the
height of $\La_N$ grows at most linearly.
It follows that the size of the largest antichain must also grow superpolynomially.
The data in Table \ref{ta3} suggest the following conjecture.

\paragraph{Conjecture.}
{\em There is an infinite sequence of values of $N$ for which the number
of dual atoms grows superpolynomially in $N$.}

In fact it seems likely that if $N =2^{n}$, a lower bound of the form
in \eqn{Eq5} (possibly with a different constant) applies to the
logarithm of the number of dual atoms,
and that for some sequence of values of $N$ a lower bound similar
to the upper bound on the right-hand side of \eqn{Eq3} will hold.
However, at present these are only conjectures.

In order to construct the orthogonal arrays needed to establish the lower bound in
Theorem \ref{th2} we make use of what we call ``mixed spreads'',
generalizing the notions of ``spread'' and ``partial spread'' from projective
geometry.
Arrays that can be constructed in this way we call ``geometric''.
Many familiar examples of orthogonal arrays, for example arrays constructed from linear codes,
are geometric.
The construction is not restricted to strength 2 (and is one of the few general
constructions we know of for mixed arrays of strength greater than 2).
The construction will be described in Section 5.

In Section 6 we use this construction to determine the lattice
$\Lambda_{64}$, and in doing so we find tight arrays with parameter sets
$$(64,2^5 4^{17} 8^1), (64, 4^{14} 8^3) ,
(64, 2^5 4^{10} 8^4), (64, 4^7 8^6) ,
$$
which appear to be new.

When studying parameter sets of putative orthogonal arrays with $N$ runs,
it is convenient to be able to say that if the number of degrees of freedom
of the parameter set $(N, s_1^{k_1} s_2^{k_2} \ldots )$, that is,
\beql{Eq61a}
k_1 (s_1 -1) + k_2 (s_2 -1) + \cdots ~,
\eeq
is small compared with $N-1$, then an orthogonal array certainly exists.

To make this precise, we define the {\em threshold function}
$B(N)$ to be the maximum number $b$ such that every
parameter set (satisfying conditions (C1) to (C4))
with at most $b$ degrees of freedom is realized by an orthogonal array, but some parameter set (again
satisfying (C1) to (C4)) with $b+1$ degrees of freedom is not realized.
If every parameter set satisfying (C1) to (C4) is realized, we set $B(N) = N-1$.

Figure \ref{fg1} shows that $B(12) =6$, since there is no $OA(12 , 2^5 3^1)$, but
every parameter set with at most 6 degrees of freedom is realized.

We are not aware of any earlier investigations of $B(N)$.

\begin{theorem}\label{th4}
If $N$ is a power of a prime then
$$N^{3/4} \le B(N) ~.$$
\end{theorem}

\noindent
In words, if the number of degrees of freedom in the parameter set
does not exceed $N^{3/4}$, then an orthogonal array exists.
This is certainly weak, but is enough to establish the lower bound of Theorem \ref{th2}.
It would be nice to have more precise estimates for $B(N)$.

A final remark.
We {\em could} have considered the partially ordered set whose nodes are all the
inequivalent orthogonal arrays with $N$ runs, rather than just their parameter sets.
However, the number of nodes then becomes unmanageably large, even for small values of $N$
(furthermore, it appears that ``meet'' and
``join'' are no longer well-defined,
and so in general this partially ordered set would not be a lattice).

Consider $N=28$, for example.
Using Kimura's (1994a, 1994b) enumeration of the Hadamard
matrices of order 28, we have calculated\footnote{Using the method described on page 165 of Hedayat, Sloane and Stufken (1999).} that
there are precisely 7570 inequivalent $OA(28, 2^{27})$'s.
This would be merely a lower bound on the number of dual atoms.
On the other hand we know (see Table \ref{ta1}) that $\La_{28}$ has precisely four dual atoms,
between 47 and 55 nodes, and height 28.

\begin{figure}[htb]
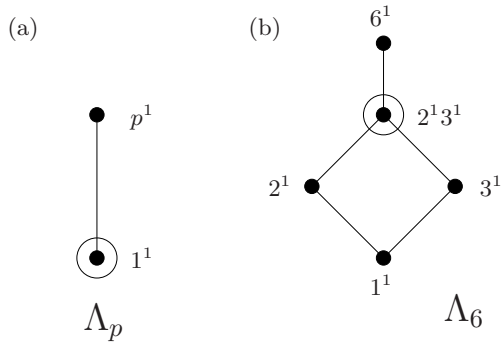

\begin{center}
\input fg2.pstex_t
\end{center}
\caption{(a)~$\Lambda_p$ and (b)~$\Lambda_6$.}
\label{fg2}
\end{figure}

\section{Examples of the lattices $\Lambda_N$ and $\Lambda'_N$}
There are a few general cases when we can describe $\La_N$ explicitly
(and for which $\La'_N$ is the same as $\La_N$).

If $N=p$ is a prime then $\La_N = \La'_N$ has two nodes, one dual
atom and height 1, as shown in Fig. \ref{fg2}(a) (dual atoms are circled).

If $N=pq$ is the product of two distinct primes, $\La_N = \La'_N$ has five nodes,
one dual atom and height 3.
$\La_6$ is shown in Fig. \ref{fg2}(b).

\begin{figure}[htb]
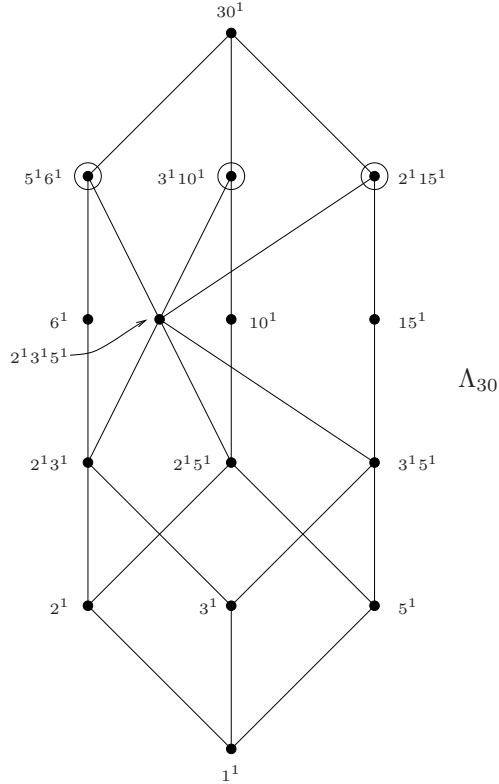

\begin{center}
\input fg3.pstex_t
\end{center}
\caption{The lattice $\Lambda_{30}$:
there are $\beta_4 = 15$ nodes, three dual atoms (circled) and the height is 5.}
\label{fg3}
\end{figure}

More generally, if $N$ is the product of $u \ge 2$ distinct primes,
it is not difficult to show that $\La_N = \La'_N$ has $2^{u-1} -1$ dual atoms,
height $2u-1$, and $\beta_{u+1}$ nodes, where
$$\{ \beta_1, \beta_2, \beta_3, \ldots \} = \{
1,2,5,15,52, \ldots \}
$$
are the Bell numbers (see Sequence \htmladdnormallink{A110}{http://www.research.att.com/cgi-bin/access.cgi/as/njas/sequences/eisA.cgi?Anum=000110} of Sloane, 1999).
Figure \ref{fg3} shows $\La_{30}$,
illustrating the case $u=3$.

\begin{figure}[htb]
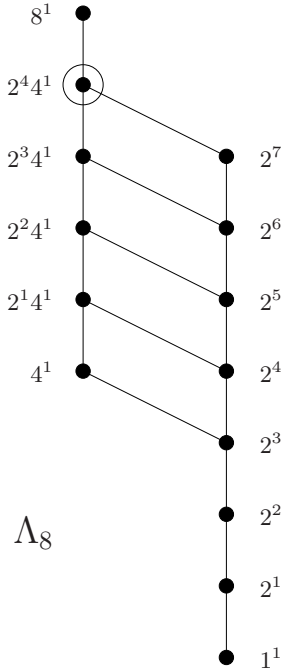

\begin{center}
\input fg4.pstex_t
\end{center}
\caption{The lattice $\Lambda_8$: 14 nodes, one dual atom, height 9.}
\label{fg4}
\end{figure}

If $N=p^2$ is the square of a prime, $\La_N = \La'_N$ has a single dual
atom, $OA(p^2, p^{p+1} )$, has height $p+2$ and contains $p+3$ nodes.
If $N=p^3$ ($p$ prime), $\La_N = \La'_N$ also has a single dual atom,
$OA(p^3, (p)^{p^2} (p^2)^1 )$, has height $p^2 + p+3$ and
contains
$2p^2 + p+4$ nodes.
$\La_8$ is shown in Fig. \ref{fg4}.
If $N= p^4$ ($p$ prime), $\La_N = \La'_N$ has two dual atoms,
$OA(p^4, (p^2)^{p^2 +1} )$ and
$OA(p^4 , (p)^{p^3} (p^3)^1)$, has height
$p^3 + 2p^2 + p+3$ and contains
$$\frac{1}{2} (p^5 + p^4 + 5p^3 + 5p^2 + 2p + 10 )$$
nodes.

For all values of $N$ mentioned so far in this section, the threshold function $B(N) = N-1$.

If $N$ is not of one of the above forms then it seems necessary to consider each case
individually.
Table \ref{ta1} summarizes the properties of $\La_N$ for some small values of $N$.
Here $A(N)$ denotes the number of dual atoms in $\La_N$.
Most of the entries in this table can be deduced from the tables in Chapter 12 
of Hedayat, Sloane and Stufken (1999),
except for $N=32$ and 64, which are discussed in Section 6.
\begin{table}[htb]
\caption{For the lattice $\La_N$ of parameter sets of orthogonal arrays with $N$ runs,
the table gives the number of dual atoms $A(N)$, the height
$ht(N)$, the total number of nodes $T(N)$ and the threshold function $B(N)$.}
$$
\begin{array}{rrrrr}
\multicolumn{1}{c}{N} &
\multicolumn{1}{c}{A(N)} &
\multicolumn{1}{c}{ht(N)} &
\multicolumn{1}{c}{T(N)} &
\multicolumn{1}{c}{B(N)} \\ \hline
1~~ & 0\quad & 0\quad & 1\quad & 0\quad \\
2~~ & 1\quad & 1\quad & 2\quad & 1\quad \\
3~~ & 1\quad & 1\quad & 2\quad & 2\quad \\
4~~ & 1\quad & 4\quad & 5\quad & 3\quad \\
5~~ & 1\quad & 1\quad & 2\quad & 4\quad \\
6~~ & 1\quad & 3\quad & 5\quad & 5\quad \\
7~~ & 1\quad & 1\quad & 2\quad & 6\quad \\
8~~ & 1\quad & 9\quad & 14\quad & 7\quad \\
9~~ & 1\quad & 5\quad & 6\quad & 8\quad \\
10~~ & 1\quad & 3\quad & 5\quad & 9\quad \\
11~~ & 1\quad & 1\quad & 2\quad & 10\quad \\
12~~ & 4\quad & 12\quad & 23\quad & 6\quad \\
13~~ & 1\quad & 1\quad & 2\quad & 12\quad \\
14~~ & 1\quad & 3\quad & 5\quad & 13\quad \\
15~~ & 1\quad & 3\quad & 5\quad & 14\quad \\
16~~ & 2\quad & 21\quad & 61\quad & 15\quad \\
17~~ & 1\quad & 1\quad & 2\quad & 16\quad \\
18~~ & 2\quad & 10\quad & 26\quad & 15\quad \\
\end{array}
$$
\label{ta1}
\end{table}

\begin{table}[htb]
$$
\begin{array}{ccccc}
\multicolumn{5}{c}{\mbox{Table 1 (cont.)}} \\ [+.2in]
\multicolumn{1}{c}{N} &
\multicolumn{1}{c}{A(N)} &
\multicolumn{1}{c}{ht(N)} &
\multicolumn{1}{c}{T(N)} &
\multicolumn{1}{c}{B(N)} \\ \hline
19 & 1 & 1 & 2 & 18 \\
20 & 4 & 20 & 35 & 11 \\
21 & 1 & 3 & 5 & 20 \\
22 & 1 & 3 & 5 & 21 \\
23 & 1 & 1 & 2 & 22 \\
24 & 4-7 & 25 & 119-133 & 18-22 \\
25 & 1 & 7 & 8 & 24 \\
26 & 1 & 3 & 5 & 25 \\
27 & 1 & 15 & 25 & 26 \\
28 & 4 & 28 & 47-55 & 15 \\
29 & 1 & 1 & 2 & 28 \\
30 & 3 & 5 & 15 & 29 \\
31 & 1 & 1 & 2 & 30 \\
32 & 2 & 42 & 320 & 29 \\
33 & 1 & 3 & 5 & 32 \\
34 & 1 & 3 & 5 & 33 \\
35 & 1 & 3 & 5 & 34 \\
\ldots & \ldots & \ldots & \ldots & \ldots \\
64 & 7 & 86 & 3037 & 57 \\
\ldots & \ldots & \ldots & \ldots & \ldots 
\end{array}
$$
\end{table}

$N=24$ is the smallest case when we do not know the complete lattice
$\La_N$.
In $\La_{24}$ the maximal value of $k$ that occurs in each of the parameter
sets $2^k 3^1$, $2^k 3^1 4^1$,
$2^k 4^1 6^1$ and $2^k 6^1$ is presently unknown.
For $2^k 3^1$, for example, it is known only that an
$OA(24, 2^{16} 3^1 )$ exists, and an $OA(24, 2^{21} 3^1)$ is impossible by the linear programming bound.
The number of dual atoms is in the range 4 to 7.
It is possible to show that the height of $\La_{24}$ is 25, however:
no chain can be longer than
$$
\mbox{$24^1$ --- $2^{20} 4^1$ --- $2^{23}$ --- $2^{22}$ --- $2^{21}$ --- $\cdots$ --- $2^1$ --- $1^1$} \,.
$$

We also do not know $\La_N$ for $N=28$, 36, $\ldots$.

The four sequences in Table \ref{ta1} are Sequences
\htmladdnormallink{A39927}{http://www.research.att.com/cgi-bin/access.cgi/as/njas/sequences/eisA.cgi?Anum=039927},
\htmladdnormallink{A39930}{http://www.research.att.com/cgi-bin/access.cgi/as/njas/sequences/eisA.cgi?Anum=039930},
\htmladdnormallink{A39931}{http://www.research.att.com/cgi-bin/access.cgi/as/njas/sequences/eisA.cgi?Anum=039931}
and
\htmladdnormallink{A48893}{http://www.research.att.com/cgi-bin/access.cgi/as/njas/sequences/eisA.cgi?Anum=048893}
of Sloane (1999).
The entries in that database will be updated as further values are determined.

We end this section with a remark about the nature of $\La_N$ as a lattice.
Since not all maximal chains between two nodes need have the same length (see Fig. \ref{fg1}),
$\La_N$ does not in general satisfy the Jordan-Dedekind chain condition (cf. Welsh, 1976).
It follows that $\La_N$ is in general not distributive, not semimodular, 
nor is it the lattice of a matroid.

\section{The maximum height of $\La_N$}
In this section we give the proof of Theorem \ref{th1}.

Let $\sigma$ denote a specification $s_1^{k_1} s_2^{k_2} \ldots$ of factors at various levels,
leaving the number of runs unspecified.
Given $\sigma$, there is a smallest number of runs,
$N_0$ say, for which an $OA(N_0, \sigma )$ exists.
Let $h$ be the height of the parameter set $(N_0, \sigma )$ in $\La_{N_0}$.
Then if the parameter set $(N, \sigma )$ occurs in any other lattice
$\La_N$, it also has height $h$.
(E.g. the specification $\sigma = 6^1$ has height 3 in each
of Figs. 1, 2(b) and 3.)
We may therefore
define $ht(\sigma)$ to be $h$, independently of the number of runs.

The height of $\La_N$, $ht(N)$, as defined in Section 1 coincides with
$ht(\sigma )$ when $\sigma = N^1$.
The height function has the following additivity property.
\begin{lemma}\label{le1}
\beql{Eq10}
ht(s_1^{k_1} s_2^{k_2} \ldots ) = k_1 ht( s_1) + k_2 ht(s_2) + \cdots ~.
\eeq
\end{lemma}

\paragraph{Proof.}
If there is a single factor on the left-hand side, say $s_1 =s$,
$k_1 =1$, then \eqn{Eq10} just says that $ht(s^1 ) = ht(s)$, repeating the assertion made above.
Otherwise, more than one factor occurs in the specification $\sigma$ on the left-hand side of \eqn{Eq10}.
Suppose the parameter set $(N, \sigma )$ occurs as a node in some $\La_N$.
The portion of $\La_N$ consisting of this node and all nodes
dominated by it is the product in an obvious sense of $k_1$
copies of $\La_{s_1}$, $k_2$ copies of $\La_{s_2}$, etc.
The height of $(N, \sigma )$ is the sum of the heights of all these
sublattices, and \eqn{Eq10} follows.~~~$\bsq$

Lemma \ref{le1} reduces the calculation of $ht (s_1^{k_1} s_2^{k_2} \ldots )$ to the calculation of the values of $ht(s_1)$, $ht(s_2) , \ldots$.
To determine the latter we must consider exactly which orthogonal arrays exist with a specified number of runs.
For this we can apply the following lemma.

\begin{lemma}\label{le2}
\beql{Eq11}
ht(N) = 1+ \max \sum_i k_i ht(s_i) ~,
\eeq
where the maximum is taken over all parameter sets 
$(N, s_1^{k_1} s_2^{k_2} \ldots ) \neq (N,N^1)$ for which an orthogonal array exists.
\end{lemma}

\paragraph{Proof.}
The height of $\La_N$ is one more than the maximal height among the dual atoms.
\eqn{Eq11} follows by applying Lemma \ref{le1} to the parameter set of such a dual atom.~~~$\bsq$

We can now use linear programming to obtain an upper bound on $ht(N)$, by maximizing
\beql{Eq12}
1+ \sum_i k_i ht (s_i)
\eeq
over all choices of $s_1$, $k_1$, $s_2$, $k_2 , \ldots$ that satisfy (C1) to (C4).

We first consider the case when $N=2^n$ for some $n$.

The case $N=64$ will illustrate the method.
If there is a factor $32^1$ then linear programming shows that \eqn{Eq12}
is maximized by $2^{32} \, 32^1$,
giving height 75.
If there is a factor $16^1$ then there is a unique parameter set that maximizes
\eqn{Eq12}, $4^{16} 16^1$, giving height 86.
Otherwise, if only factors $2^{k_1} 4^{k_2} 8^{k_3}$ occur, the height does not
exceed 85.
Since an $OA(64, 4^{16} 16^1 )$ exists, we conclude that $ht( 64) = 86$.
We will return to the case $N=64$ in Section 6.

In this way we obtain the values of $ht(N)$, $N=2^n$ shown in Table \ref{ta2}.
\begin{table}[htb]
\caption{Height of $\La_N$ for $N=2^n$.}
$$
\begin{array}{cccccc}
n & 1 & 2 & 3 & 4 & 5 \\
N=2^n & 2 & 4 & 8 & 16 & 32 \\
ht(N) & 1 & 4 & 9 & 21 & 42 \\
ht(N) /(N-1) & 1 & 1.3333 \ldots & 1.2857 \ldots & 1.4 & 1.3548 \ldots \\ [+.5in]
n & 6 & 7 & 8 & 9 & 10 \\
N=2^n & 64 & 128 & 256 & 512 & 1024 \\
ht(N) & 86 & 171 & 358 & 715 & 1431 \\
ht(N) / (N-1) & 1.3650 \ldots & 1.3465 \ldots & 1.4039 \ldots & 1.3992 \ldots & 1.3988 \ldots
\end{array}
$$
\label{ta2}
\end{table}

Consider the general problem of maximizing \eqn{Eq12} for $N=2^n$.
Comparing \eqn{EqR} and \eqn{Eq12}, we see that an $s$-level
factor contributes $ht(s)$ to the height but uses up $s-1$ degrees of freedom.
If we ignore the constraints of integrality then the value of the expression
in \eqn{Eq12} would be maximized
by a term $s^k$ where $s$ is chosen to maximize $ht(s) / (s-1)$.
This suggests that we should investigate this quantity in order to prove
Theorem \ref{th1}.
The data in
Table \ref{ta2} suggest that the ratio $ht(s) / (s-1)$ is maximized if $s$ is of the form
$2^{2^i}$ and is as large as possible.
That this is indeed so is established by the next three lemmas.

\begin{lemma}\label{L200}
Given positive real numbers $\af_r$, $\beta_r$ $(r=1,\ldots, m )$, $\gamma$,
the maximal value of
$$
\sum_{r=1}^m \alpha_r x_r
$$
subject to the constraints
$$
\begin{array}{c}
\displaystyle\sum_{r=1}^m \beta_r x_r = \gamma ~, \\[+.2in]
x_r \ge 0 \quad (r=1, \ldots, m )
\end{array}
$$
is
$$\gamma \max_{r=1 \ldots m} \frac{\af_r}{\beta_r} ~.
$$
\end{lemma}

We omit the straightforward proof.
\begin{lemma}\label{L201}
Given positive real numbers $\alpha_r$, $\beta_r$ $(r=1, \ldots, m+n )$,
$\gamma$, $\gamma'$ with $\gamma' \le \gamma$, the maximal
value of
\beql{Eq202}
\sum_{r=1}^{m+n} \alpha_r x_r
\eeq
subject to the constraints
\beql{Eq203A}
\sum_{r=1}^{m+n} \beta_r x_r = \gamma ~, 
\eeq
\beql{Eq203B}
\sum_{r=m+1}^{m+n} \beta_r x_r \le \gamma' ~, 
\eeq
\beql{Eq203C}
x_r \ge 0 \quad (r=1, \ldots, m+n)
\eeq
is
\beql{Eq204}
\max \left\{
\gamma \max_{r=1 \ldots m} \frac{\af_r}{\beta_r} ,~
(\gamma - \gamma' ) \max_{r=1 \ldots m}
\frac{\af_r}{\beta_r} +
\gamma' \max_{r = m+1 \ldots m+n} \frac{\af_r}{\beta_r} \right\}\,.
\eeq
The maximum is given by the first expression if and only if
$$
\max_{r=1 \ldots m} \frac{\af_r}{\beta_r} \ge
\max_{r= m+1 \ldots m+n} \frac{\af_r}{\beta_r} ~.
$$
\end{lemma}

\paragraph{Proof.}
Let $\sigma$, $0 \le \sigma \le \gamma'$, denote the value of the left-hand
side of \eqn{Eq203B}.
Then by Lemma \ref{L200} the maximal value of the sum in \eqn{Eq202} is equal to
$$(\gamma - \sigma ) \max_{r=1 \ldots m}
\frac{\af_r}{\beta_r} + \sigma
\max_{r= m+1 \ldots m+n} \frac{\af_r}{\beta_r} \,.
$$
This is a linear function of $\sigma$ and so its maximal value is taken
at one of the two endpoints, leading to \eqn{Eq204}.~~~$\bsq$

We can now give an upper bound on the height $ht(N)$ for $N=2^n$.
Let
$$\rho_n = \frac{ht(N)}{2^n -1} , \quad N=2^n ~.$$

\begin{lemma}\label{L205}
$$\rho_n < c \quad\mbox{for all}\quad n ~,$$
where $c$ is the constant
$$\sum_{i=0}^\infty \frac{1}{2^{2^i} -1} = 1.4039 \ldots$$
\end{lemma}

\paragraph{Proof.}
We use the induction hypothesis that, for $m \ge 0$,
\beql{Eq206}
\rho_{2^m} > \rho_x \quad\mbox{whenever}\quad
1 \le x < 2^{m+1}, ~x \neq 2^m \,.
\eeq
This is trivially true when $m=0$.
We first compute $\rho_{2^m}$.
From the Rao-Hamming construction, for $N=2^{2^m}$, an $OA(N, (\sqrt{N})^{\sqrt{N} +1} )$ always exists, and we find that
$$
\rho_{2^m} \ge \rho_{2^{m-1}} + \frac{1}{2^{2^m} -1} ~.
$$
On the other hand, we obtain an upper bound on $\rho_{2^m}$ from the linear program:
choose nonnegative integers $k_1, k_2, \ldots$ so as to maximize
$$
1 + \sum_{r=1}^{2^m -1} ht (2^r) k_{r}
$$
subject to the constraint
$$
\sum_{r=1}^{2^m -1} (2^r -1) k_{r} = 2^{2^m} -1 ~.
$$
From Lemmas \ref{le2} and \ref{L200} we have
\begin{eqnarray*}
\rho_{2^m} & \le & \frac{1}{2^{2^m} -1} + \max_{r=1 \ldots 2^m -1} \rho_r \\
& = & \frac{1}{2^{2^m} -1} + \rho_{2^{m-1}} ~.
\end{eqnarray*}
Since the two bounds agree,
\beql{Eq209}
\rho_{2^m} = \rho_{2^{m-1}} + \frac{1}{2^{2^m} -1} ~.
\eeq


We now complete the proof of the induction step.
For $1 \le x < 2^m$, we have
$$\rho_{2^m} > \rho_{2^{m-1}} \ge \rho_x ~,$$
as required.
Suppose $2^m < x < 2^{m+1}$.
Then $\rho_x$ is upper-bounded by the solution to the linear program:
maximize
$$\left(1 + \sum_{r=1}^{x-1} ht (2^r ) k_{r} \right) \Bigl/ (2^x -1 )$$
subject to the constraints
\begin{eqnarray*}
\sum_{r=1}^{x-1} (2^r -1) k_{r} & = & 2^x -1 ~,
\\
\sum_{r=2^m}^{x-1} (2^r -1) k_{r} & \le & 2^{x-1} -1 ~.
\end{eqnarray*}
(The second constraint is implied by the requirement that there can be
at most one factor which has more levels than the square root of the number of runs.)
By Lemma \ref{L201} and induction on $x$, we obtain
\begin{eqnarray*}
\rho_x & \le & \frac{1}{2^x -1}
( 1+ 2^{x-1} \rho_{2^{m-1}} + (2^{x-1} -1) \rho_{2^m} ) \\
& = & \rho_{2^m} - \frac{1+ 2^{x-1} - 2^{2^m}}{(2^{2^m} -1) (2^x -1)} \\
& < & \rho_{2^m} ~.
\end{eqnarray*}
This completes the induction step.

To complete the proof of Lemma \ref{L205}, by the induction hypothesis it suffices to prove that
$\rho_{2^m} < c$ for all $m$.
But from \eqn{Eq209} it follows that
$$\rho_{2^m} = \sum_{i=0}^m \frac{1}{2^{2^i} -1} < c ~,$$
and that $\rho_{2^m} \to c$ as $m \to \infty$.~~~\bsq

If $N$ is not a power of 2 then similar arguments show that
%
%
the height is
(considerably) less than $cN$.
This completes the proof of Theorem \ref{th1}.

\section{Upper bounds on the number of parameter sets}
In this section we establish the upper bounds in
Theorems \ref{th2} and \ref{th3}.
We will bound $T'(N)$, the number of nodes in $\La'_N$.
Since $\La_N$ is a sublattice of $\La'_N$, this is also an upper bound
on the number of nodes in $\La_N$.
Suppose first that $N=2^{2r}$.

We start by considering parameter sets $(N, 2^{k_1} 4^{k_2} 8^{k_3} \ldots (2^r)^{k_r} )$, containing no level exceeding $\sqrt{N}$.
From \eqn{EqR},
\beql{Eq19}
N-1 \ge k_1 + 3k_2 + 7k_3 + \cdots +
(2^r -1) {k_r} ~.
\eeq
Let $\gamma$ be the number of nonnegative integer solutions
$(k_1, k_2, \ldots, k_N )$ to this inequality.
Then $\gamma / (N-1)^r$ is the Riemann sum approximating the volume
of the simplex bounded by the hyperplanes
\begin{eqnarray*}
&& 1 \ge x_1 + 3x_2 + 7x_3 + \cdots + (2^r -1) x_r ~, \\
&& x_1 \ge 0 , ~ x_2 \ge 0 , \ldots, x_r \ge 0
\end{eqnarray*}
in $\RR^r$.
For large $N$ this yields
\beql{Eq20}
\gamma = \frac{N^r}{r! \prod\limits_{i=1}^r (2^i -1)}
(1+o(1)) ~.
\eeq
The product in the denominator approaches $c_3 2^{r(r+1)/2}$, as $r \to \infty$, where $c_3 = 0.2887 \ldots$.

Now suppose the
parameter set contains a factor at $2^i$ levels, where $r+1 \le i \le 2r-1$.
There can be at most one such factor, and the number of such parameter
sets in each case is at most $\gamma$.
The total number of parameter sets is therefore at most $r \gamma$, and setting $r= \frac{1}{2} \log_2 N$ we find that 
$$ \log_2 (r \gamma ) \le \frac{3}{8} (\log_2 N)^2 (1+o(1)) ~.$$
This establishes the upper bound in Theorem \ref{th2} for $N=2^{2r}$. It also implies 
the upper bound for $N=2^{2r+1}$, after noting that $T'(N) \le T'(2N)$.

We now give a sketch of the proof of Theorem \ref{th3}, omitting many tedious
details.
To simplify the analysis we will 
neglect terms on the right-hand side of (\ref{EqR}) that correspond to factors 
with a level greater than $\sqrt{N}$.
Suppose first that $N$ is a large number of the form $2^{2a_1} 3^{2a_2}$.
Pretending for the moment that $a_1$ and $a_2$ are allowed to be real numbers, not just integers, we may consider what choice of $a_1$ and $a_2$ maximizes the number of solutions to \eqn{EqR} for a given value of $N$.
The number of terms on the right-hand side of \eqn{EqR} is now $(a_1 + 1)(a_2 + 1) - 1$.
The arguments used to establish the upper bound of Theorem \ref{th2} show that
the number of solutions to \eqn{EqR}
is maximized if $2^{2a_1}$ is approximately equal to $3^{2a_2}$.

Now suppose that $N$ is of the form
\beql{Eq7A}
p_1^{2a_1} ~ p_2^{2a_2} \cdots p_m^{2a_m} ~,
\eeq
where $p_1 =2$, $p_2=3, \ldots$ are the first $m$ primes.
We find that
the number of solutions to \eqn{EqR} is maximized when the numbers
$p_i^{2a_i}$ are all approximately equal, and we will therefore assume that $p_i^{2a_i} = N^{\frac{1}{m} (1+o(1))}$, i.e. that
$$a_i = \frac{1}{2m} \frac{\ln N}{\ln p_i} (1+o(1)) , \quad i =1, \ldots, m ~.$$

The Rao bound contains a term for every possible level
$$s = p_1^{i_1} ~ p_2^{i_2} \cdots p_m^{i_m}$$
in which $0 \le i_{\nu} \le a_\nu$, $1 \le \nu \le m$, where not
all the $i_\nu$ are equal to 0.
The number of such terms is
\beql{Eq30}
\delta  :=  (a_1 + 1) (a_2 + 1) \cdots (a_m + 1) -1 =
\frac{1}{(2m)^m} \frac{(\ln N)^m}{\prod\limits_{j=1}^m \ln p_j}
(1+ o(1)) ~.
\eeq
The product of the coefficients of all the terms on the right-hand side of the Rao bound is
$$\zeta := p_1^{a_1^2 a_2 \ldots a_m /2}
p_2^{a_1 a_2^2 a_3 \ldots a_m /2} \cdots
p_m^{a_1 \ldots a_{m-1} a_m^2/2} (1+ o(1)) ~.
$$
This implies
$$\ln \zeta = \frac{\delta}{4} \ln N (1+o(1)) ~.$$
Again using $\gamma$ to denote the number of solutions to the Rao inequality, we have
$$\gamma = \frac{N^\delta}{\delta ! \zeta} (1+ o(1)) ~,$$
hence
\begin{eqnarray*}
\ln \gamma & = & \frac{3}{4} \frac{(\ln N)^{m+1}}{(2m)^m (\ln m)^{m (1+o(1))}} -
\frac{1}{(2m)^m} \frac{(\ln N)^m}{(\ln m)^{m(1+o(1))}} (m \ln \ln N - m \ln m) \\
&&~~~~+~ \mbox{smaller terms} ~.
\end{eqnarray*}
This expression is maximized if we take
$$m = \frac{1}{2e} \frac{\ln N}{\ln \ln N} (1+o(1)) ~,$$
and then we find that the leading term in the expression for $\ln \gamma$ is
$$\frac{3}{4} \ln N N^{\frac{1}{2e} \frac{1}{\ln \ln N}} ~.$$
We conclude that
$$\ln \ln \gamma \le \frac{1}{2e} \frac{\ln N}{\ln \ln N} ~,$$
which establishes Theorem \ref{th3}.

\section{Geometric orthogonal arrays}
We consider subspaces $V$ of the vector space $GF(q)^n$ over $GF(q)$, where $q$ is a power of a prime.
By the
dimension of $V$, $\dim V$, we mean the vector space dimension over $GF(q)$ (rather than the projective dimension, which is one less).
The following
notion was suggested by the notions of spread and partial spread in projective geometry
(cf. Thas, 1995).

\paragraph{Definition.}
A {\em mixed spread of strength $t$} is a collection
$\sV = \{V_1, V_2 , \ldots, V_k \}$ of subspaces of $GF(q)^n$ such that for all choices
of $\tau \le t$ indices $i_1, i_2, \ldots, i_\tau$ 
(with $1 \le i_1 < i_2 < \cdots < i_\tau \le k$) the dimension
of the span of $V_{i_1}, \ldots, V_{i_\tau}$ is equal to $\dim V_{i_1} + \cdots + \dim V_{i_\tau}$.

An equivalent condition is that the span of $V_{i_1}, \ldots, V_{i_\tau}$ is
the direct sum $V_{i_1} \oplus \cdots \oplus V_{i_\tau}$
for all choices of $\tau \le t$ indices $i_1, \ldots, i_\tau$ 
with $1 \le i_1 < \cdots < i_\tau \le k$.

Any collection $\sV$ of subspaces has strength 1.
$\sV$ has strength 2 if and only if every pair $V_i$, $V_j \in \sV$, $i \neq j$, intersect just in the zero vector.
$\sV$ has strength 3 if and only if it has strength 2 and for any triple of distinct
subspaces each one meets the span of the other two just in the zero vector.

If $V$ is a $d$-dimensional subspace of $GF(q)^n$ we denote by $V^\ast$ the dual space, the
space of linear functionals on $V$
(see for example Hoffman and Kunze, 1961),
and we fix a labeling $f_0 , f_1, \ldots, f_{q^d -1}$ for the elements of $V^\ast$.

Given a mixed spread of strength $t$, $\sV= \{ V_1, V_2, \ldots, V_k \}$, where the $V_i$ are 
subspaces of $GF(q)^n$, we obtain an orthogonal array $OA(\sV)$ with $q^n$ runs as follows.
The columns of the array are labeled $V_1, V_2, \ldots, V_k$ and the rows are labeled
by the linear functionals $f \in (GF(q)^n )^\ast$.
If $f$ restricted to $V_i$, $f |_{V_i}$, is the $j$th linear functional in $V_i^\ast$,
the $(f, V_i)$-th entry in the array is $j$.
The symbols in column $i$ are therefore taken from $\{0,1, \ldots, q^{\dim V_i} -1 \}$.

We will say that an orthogonal array constructed in this way is {\em geometric}.

\begin{theorem}\label{P1}
The orthogonal array $OA(\sV)$ has strength $t$ if and only if the mixed spread $\sV$ has strength $t$.
\end{theorem}

\paragraph{Proof.}
Suppose $OA(\sV)$ has strength $t$.
Consider for example the first $t$ columns.
In the projection of the array onto these columns we see
$$\prod_{i=1}^t q^{\dim V_i} = q^{\sum\limits_{i=1}^t \dim V_i}$$
different $t$-tuples of symbols.
Since these depend only on the restrictions of the $f \in (GF(q)^n)^\ast$ to the span of
$V_1, \ldots, V_t$, the dimension of that space must be at least the sum of the dimensions of $V_1, \ldots, V_t$, and clearly it cannot have a higher dimension.
So $\sV$ is a mixed spread of strength $t$.

Conversely, suppose $\sV$ is a mixed spread of strength $t$.  We can write
$$GF(q)^n = V_{1} \oplus \cdots \oplus V_{t} \oplus X$$
where $X$ is the complementary space to the $V_i$.
Since the dual of a direct sum is canonically isomorphic to the direct sum of
the duals, we immediately find that as we run through the
linear functionals on $GF(q)^n$, every tuple
$(f|_{V_1} , \ldots f|_{V_t} , f|_X )$ 
of restrictions occurs precisely once.
Ignoring the last component, we see that every tuple
$(f|_{V_1} , \ldots f|_{V_t} )$ occurs precisely $|X|$ times.
Hence $OA(\sV)$ has strength $t$.~~~$\bsq$

\begin{lemma}\label{P5}
Any geometric array of strength $2$ can always be extended to a tight array (i.e. one meeting the Rao bound) by adding $q$-level factors.
\end{lemma}

\paragraph{Proof.}
We simply group any unused points into 1-dimensional subspaces.~~~$\bsq$

\paragraph{Examples.} (i)
The 1-dimensional subspaces of $GF(q)^n$ form a mixed spread of strength 2.
The corresponding array is the familiar
$$OA(q^n , q^k ) , \quad k = (q^n -1) /(q-1) ~,$$
of the Rao-Hamming construction.

(ii) More generally, a classical $a$-spread in $PG(b,q)$ is a mixed spread of strength 2 in our sense.
This is a set of subspaces of $PG(b,q)$ of projective dimension $a$ which
partitions $PG(b,q)$ (Thas, 1995), and exists if and only if $a+1$ divides $b+1$.
From Theorem \ref{P1} we obtain an
$$OA(q^{b+1}, (q^{a+1} )^k ), \quad
k = (q^{b+1} -1)/(q^{a+1} -1) \,,
$$
which of course is also given by the Rao-Hamming construction.

We could also have obtained example (ii) directly from example (i), by remarking that a mixed
spread of strength $t$ over $GF(q)$, $q =p^\beta$, is also a mixed spread
of strength $t$ over $GF(q')$, $q' = p^\af$, if $q'$ divides $q$.
The dimensions of the subspaces are multiplied by $\beta / \af$.

(iii) Provided $a \ge b/2$, there exists a mixed spread of strength 2
in $GF(q)^b$ consisting of a single subspace $GF(q)^a$ and a partitioning of the remaining points into $q^a$ subspaces $GF(q)^{b-a}$.
This can be proved directly, or alternatively is equivalent to Lemma 2.1
of Eisfeld, Storme and Sziklai (1999).
From Theorem \ref{P1} we obtain a geometric
$$OA(q^b, (q^{b-a} )^{q^a} (q^a)^1 )$$
whenever $a \ge b/2$.
Orthogonal arrays with these parameters were already known from the difference scheme
construction (Hedayat, Sloane and Stufken, 1999, Example 9.19), but the present
construction also shows that they are geometric.

(iv) The classical ``partial $a$-spread'' constructed in Lemma 2.2 of Eisfeld, Storme and Sziklai (1999)
translates in our language into a mixed spread of strength 2 consisting of $k$ $b$-dimensional
subspaces $(b \ge 2)$ of $GF(q)^n$, where $n = ib+r$, $0 \le r< b$, and
$$k =q^r \frac{q^{ib} -1}{q^b -1} - q^r +1 ~.$$
This produces a geometric $OA(q^n, (q^b)^k)$
(again arrays with these parameters were known from
the difference scheme construction),
which by Lemma \ref{P5} can be extended to a tight
\beql{Eq80}
OA(q^n, (q^1)^l (q^b)^k ) ~,
\eeq
where $l=q^b (q^r -1)/(q-1)$.
The orthogonal arrays constructed by Wu (1989) are a special case of \eqn{Eq80},
but in general these arrays may be new.

(v) Generalizing examples (i) and (ii), any orthogonal array formed from the codewords of
a projective linear code (one for which the columns of a generator matrix are nonzero
and projectively distinct) is geometric.

(vi) The $OA(256, 2^{16} )$ of strength 5 formed from the Nordstrom-Robinson code (see Hedayat, Sloane and Stufken, 1999, Section 5.10) is not geometric, and no geometric
$OA(256, 2^{16} )$ of strength 5 exists.

We shall see other examples in Section 6.

\paragraph{Remarks.}
An unmixed geometric orthogonal array is always linear, in the sense of Hedayat, Sloane and Stufken (1999), Chapter 3.
In general a mixed geometric orthogonal array is additive but not necessarily linear\footnote{For the distinction between additive and linear sets in the
context of coding theory see Calderbank et~al. (1998).}
over each of the fields involved.

If the strength is 2, the number of degrees of freedom in the parameter set for $OA(\sV)$ is equal to the total number of nonzero points in all the subspaces $V_i$.

Finally, the following is a recipe for constructing the orthogonal
array from a mixed spread 
$\sV = \{V_1, V_2 , \ldots, V_k \}$ of subspaces of $GF(q)^n$
in the case when $q$ is a prime.
Let $v_{1}^{(i)}, \ldots, v_{d_i}^{(i)}$ be a basis for $V_i$,
where $d_i = \dim V_i$, $1 \le i \le k$.
Let $w_0, \ldots, w_{q^n - 1}$ be the vectors of $GF(q)^n$.
Then the $i$th entry of the $j$th row of the orthogonal array,
for $1 \le i \le k$, $0 \le j \le q^n - 1$, 
is the number 
$$ \sum_{r=1}^{d_i} w_j \cdot v_r^{(i)} ~ q^{r - 1} ~.$$
(This is a number in the range $\{0, \ldots, q^{d_i}-1\}.)$

\section{If the number of runs is a power of 2}
In this section we consider the case $N=2^n$, $n=1,2, \ldots$.
We have already discussed $ht(N)$ in Section 3 (see Table \ref{ta2}).
With the assistance of Michele Colgan, we used a computer to determine the number of dual atoms $A' (N)$ and the total number of nodes $T' (N)$ in the idealized lattice $\La'_N$ for $n \le 9$.
The results are shown in the second and third columns of Table \ref{ta3}.
Note in particular the extremely rapid growth from $N=256$ to $N=512$.
We regard this as convincing evidence that when $N=2^n$, $A'(N)$ (and therefore
presumably $A(N)$) grows faster than any polynomial in $N$.
\begin{table}[htb]
\caption{Dual atoms and total number of nodes in idealized lattice $\La'_N$
$(A'(N)$ and $T' (N))$ and in lattice $\La_N$ $(A(N)$ and $T(N))$.}
$$
\begin{array}{rrrrr}
\multicolumn{1}{c}{N} &
\multicolumn{1}{c}{A'(N)} &
\multicolumn{1}{c}{T'(N)} &
\multicolumn{1}{c}{A(N)} &
\multicolumn{1}{c}{T(N)}
\\ [+.1in]
1 & 0 & 1 & 0 & 1 \\
2 & 1 & 2 & 1 & 2 \\
4 & 1 & 5 & 1 & 5 \\
8 & 1 & 14 & 1 & 14 \\
16 & 2 & 61 & 2 & 61 \\
32 & 3 & 322 & 2 & 320 \\
64 & 11 & 3058 & 7 & 3037 \\
128 & 21 & 33364 \\
256 & 72 & 789085 \\
512 & 144521 & 18614215
\end{array}
$$
\label{ta3}
\end{table}

As to the lattice $\La_N$ itself, for $n \le 4$ this is covered by the results in Section 2.
For $N=32$ there are precisely two parameter sets in $\La'_{32}$ which do not exist, $(32, 4^{10} )$ and $(32, 2^1 4^{10} )$.
These can be ruled out either by the linear programming bound or by the Bose-Bush bound (Hedayat, Sloane and Stufken,
1999, Theorem 2.8).
All other parameter sets in $\La'_{32}$ are realized.
It follows that $\La_{32}$ contains exactly two dual atoms, $OA(32, 2^{16} 16^1 )$ and $OA(32, 4^8 8^1 )$.

Before considering $\La_{64}$ we give a lemma that will be used to construct new arrays.

\begin{lemma}\label{P2}
Suppose $V_1$, $V_2$, $V_3$ are three $r$-dimensional subspaces of $GF(2)^{2r}$ such that
$V_i \cap V_j = \{0\}$, $i \neq j$.
Then their union can be replaced by $2^r -1$ two-dimensional subspaces,
any pair of which meet just in the zero vector.
\end{lemma}

{\em Proof.}
Since $V_1 \cap V_2 = \{ 0\}$, $V_1$ and $V_2$ span the space
$GF(2)^{2r}$.
Let $\pi_1$, $\pi_2$ be the associated projection maps from $GF(2)^{2r}$ to $V_1$, $V_2$ respectively.
Then $i_1 = \pi_1 \bigl|_{V_3} : V_3 \to V_1$ and $i_2 = \pi_2 \bigl|_{V_{3}} : V_3 \to V_2$ are both
isomorphisms.
It follows that $V_3$ is the set
$$\{ v+ i(v) : v ~\mbox{in}~
V_1 \} , ~~i=i_2 i_1^{-1} ~.
$$
But then we need simply take the planes $\{0, v, i(v), v+i(v) \}$ for $v \in V_1$ to establish
the lemma.~~~$\bsq$

The lemma implies that if a geometric $OA(2^{2r}, \ldots )$ exists then so does the array
obtained by replacing $(2^r)^3$ in the parameter set by $4^k$, $k=2^r -1$.
In particular, in a geometric $OA(64, \ldots )$ we can replace $8^3$ by $4^7$.

\begin{theorem}\label{P3}
The lattice $\La_{64}$ contains precisely seven dual atoms, with parameter sets
\begin{eqnarray}\label{Eq69}
&& (64, 2^5 4^{17} 8^1), (64, 4^{14} 8^3 ) , (64, 2^5 4^{10} 8^4) , \nonumber \\
&& (64, 4^7 8^6 ), (64, 8^9), (64, 4^{16} 16^1 ) , (64, 2^{32} 32^1 ) \,.
\end{eqnarray}
A geometric orthogonal array exists for each of these parameter sets.
\end{theorem}

\paragraph{Proof.}
As an intermediate step, we use
mixed spreads of strength 2
to construct orthogonal arrays with the following parameter sets:
\begin{eqnarray}\label{Eq70}
&& (64, 4^{21} ),
(64, 2^5 4^{17} 8^1),
(64, 2^4 4^{15} 8^2),
(64, 4^{14} 8^3) , \nonumber \\
&& (64, 2^5 4^{10}  8^4) ,
(64, 2^4 4^8 8^5) , (64, 4^7 8^6) , (64, 2^8 4^2 8^7) , \nonumber \\
&& (64, 2^4 4^1 8^8),
(64, 8^9) , (64, 4^{16} 16^1),
(64, 2^{32} 32^1 ) \,.
\end{eqnarray}
On the other hand, linear programming shows that orthogonal arrays do not exist with
parameter sets
\beql{Eq71}
(64, 4^{18} 8^1), (64, 4^{16} 8^2) , (64, 4^{11} 8^4), (64, 4^9 8^5),
(64, 4^3 8^7), (64, 4^2 8^8) \,.
\eeq
We then check that every parameter set with 64 runs either dominates one of \eqn{Eq71}
(and so is not realized), or is dominated by one of \eqn{Eq70} (and is realized).
Furthermore, the parameter sets in \eqn{Eq69} dominate all of \eqn{Eq70}.

It remains to construct the arrays mentioned in \eqn{Eq69}.
The last two follow from Example (iii) of Section 5.
Also $(64, 8^9)$ comes from Example (i), and $(64, 4^7 8^6)$ and
$(64, 4^{14} 8^3)$ follow from Lemma \ref{P2}.

To construct an $OA(64, 2^5 4^{10} 8^4)$ we proceed as follows.
We begin by constructing an explicit example of an $OA( 64, 8^9)$ from Theorem \ref{P1}, by using an extended Reed-Solomon code of length 9, dimension 2
and minimal distance 8 over $GF(8)$.
This gives a decomposition of $GF(2)^6$ into 9 copies of $GF(2)^3$
meeting only in the zero vector.
These 9 subspaces are spanned by
the following nine triples of columns:
\beql{Eq72}
\begin{array}{ccccccccc}
0 & I & I & I & I & I & I & I & I \\
I & 0 & I & A & A^2 & A^3 & A^4 & A^5 & A^6
\end{array}
\eeq
where
$$0 = \left[ \begin{array}{ccc}
0 & 0 & 0 \\
0 & 0 & 0 \\
0 & 0 & 0
\end{array}
\right], \quad
I = \left[
\begin{array}{ccc}
1 & 0 & 0 \\
0 & 1 & 0 \\
0 & 0 & 1
\end{array}
\right] , \quad
A = \left[
\begin{array}{ccc}
0 & 1 & 0 \\
1 & 0 & 1 \\
0 & 1 & 1
\end{array}
\right] ,
$$
and $A^7 = I$.

We may replace the first four subspaces and the last subspace (which together
contain 35 nonzero points) by ten two-dimensional
subspaces with five single points left over.
One choice for the ten two-dimensional subspaces
is shown in Table \ref{ta4}.
\begin{table}[htb]
\caption{Ten pairs of columns each spanning a two-dimensional
subspace of $GF(2)^6$.}
$$
\begin{array}{c@{\,}c}
0 & 1 \\
0 & 0 \\
0 & 1 \\
0 & 0 \\
1 & 0 \\
1 & 1 \\
\end{array}
\qquad
\begin{array}{c@{\,}c}
0 & 1 \\
0 & 1 \\
0 & 0 \\
1 & 1 \\
0 & 1 \\
1 & 0 \\
\end{array}
\qquad
\begin{array}{c@{\,}c}
0 & 0 \\
0 & 1 \\
1 & 0 \\
0 & 0 \\
1 & 1 \\
1 & 0 \\
\end{array}
\qquad
\begin{array}{c@{\,}c}
0 & 0 \\
0 & 0 \\
0 & 1 \\
0 & 0 \\
0 & 0 \\
1 & 0 \\
\end{array}
\qquad
\begin{array}{c@{\,}c}
0 & 1 \\ 0 & 1 \\ 0 & 0 \\ 1 & 0 \\ 1 & 0 \\ 1 & 0 \\
\end{array}
\qquad
\begin{array}{c@{\,}c}
0 & 1 \\ 1 & 1 \\ 1 & 1 \\ 0 & 1 \\ 0 & 1 \\ 0 & 1 \\
\end{array}
\qquad
\begin{array}{c@{\,}c}
0 & 0 \\ 0 & 1 \\ 0 & 0 \\ 1 & 0 \\ 0 & 0 \\ 0 & 0 \\
\end{array}
\qquad
\begin{array}{c@{\,}c}
0 & 1 \\ 0 & 0 \\ 0 & 0 \\ 1 & 0 \\ 1 & 1 \\ 0 & 0 \\
\end{array}
\qquad
\begin{array}{c@{\,}c}
1 & 0 \\ 0 & 1 \\ 0 & 1 \\ 0 & 1 \\ 0 & 1 \\ 0 & 0 \\
\end{array}
\qquad
\begin{array}{c@{\,}c}
0 & 1 \\ 1 & 1 \\ 0 & 1 \\ 1 & 0 \\ 0 & 0 \\ 1 & 0
\end{array}
$$
\label{ta4}
\end{table}
This gives a mixed spread of strength 2 consisting of four 3-dimensional subspaces,
ten two-dimensional subspace and five points,
and so by Theorem \ref{P1} corresponds to an $OA(64, 2^5 4^{10} 8^4 )$.
Finally, Lemma \ref{P2} produces an $OA(64, 2^5 4^{17} 8^1 )$.~~~$\bsq$

\paragraph{Remark.}
Geometric orthogonal arrays with parameter sets of the form
$(64, \ldots 8^k \ldots )$ involve selecting $k$ disjoint (except for the zero vector) copies of $GF(2)^3$ inside $GF(2)^6$.
It is simpler to work projectively, and then we must choose $k$ disjoint copies of $PG(2,2)$ inside $PG(5,2)$.
Equation \eqn{Eq72} then gives a decomposition of $PG(5,2)$ into nine copies of $PG(2,2)$.

With the help of Magma (Bosma and Cannon, 1995;
Bosma, Cannon and Mathews, 1994;
Bosma, Cannon and Playoust, 1997),
we showed that if $1 \le k \le 4$ there is a unique way to choose $k$ disjoint
$PG(2,2)$'s in $PG(5,2)$, and these are equivalent to a subset of \eqn{Eq72}.
For $k=5$, there are precisely two ways, one of which is equivalent
to a subset of \eqn{Eq72}
while the other contains no $PG(2,2)$ in its complement, and so cannot be extended to $k=6$.
For $k=6, \ldots, 9$, there is again a unique way to choose $k$ disjoint planes.
In particular the decomposition into nine planes shown in \eqn{Eq72} is also unique.

An example of a maximal set of five $PG(2,2)$'s in $PG(5,2)$ is shown in Table \ref{ta5}.
This corresponds to a geometric $OA(64, 8^5)$ that cannot be extended to a geometric $OA(64, 8^6)$.
It would be interesting to determine if it can be extended to {\em any} $OA(64, 8^6)$.
\begin{table}[htb]
\caption{A set of five disjoint $PG(2,2)$'s in $PG(5,2)$ that is not contained in a set of six.  Each triple of columns spans one of the subspaces.}

$$
\begin{array}{c@{\,}c@{\,}c}
1 & 0 & 0 \\
0 & 1 & 0 \\
0 & 0 & 1 \\
0 & 0 & 0 \\
0 & 0 & 0 \\
0 & 0 & 0 \\
\end{array}
\qquad
\begin{array}{c@{\,}c@{\,}c}
0 & 0 & 0 \\
0 & 0 & 0 \\
0 & 0 & 0 \\
1 & 0 & 0 \\
0 & 1 & 0 \\
0 & 0 & 1 \\
\end{array}
\qquad
\begin{array}{c@{\,}c@{\,}c}
1 & 0 & 0 \\
0 & 1 & 0 \\
0 & 0 & 1 \\
1 & 0 & 0 \\
0 & 1 & 0 \\
0 & 0 & 1 \\
\end{array}
\qquad
\begin{array}{c@{\,}c@{\,}c}
1 & 0 & 0 \\
0 & 1 & 0 \\
0 & 0 & 1 \\
0 & 1 & 0 \\
1 & 0 & 1 \\
0 & 1 & 1 \\
\end{array}
\qquad
\begin{array}{c@{\,}c@{\,}c}
1 & 0 & 0 \\
0 & 1 & 0 \\
0 & 0 & 1 \\
1 & 1 & 1 \\
1 & 1 & 0 \\
1 & 0 & 0 \\
\end{array}
$$
\label{ta5}
\end{table}


\section{The existence of orthogonal arrays with certain parameter sets}
In this section we prove Theorem \ref{th4}, the lower bound in Theorem \ref{th2}, and also give some other conditions which are sufficient to guarantee that a parameter set can be realized by an orthogonal array.

\begin{lemma}\label{le6}
Suppose $N=p^m$ is a power of a prime and $(N, s_1^{k_1} s_2^{k_2} \ldots )$
is a parameter set with $k= \Sigma_i k_i$ factors.
If $k \le p^{\lfloor (m+1)/2 \rfloor} +1$ then this parameter set is realized by a geometric orthogonal array.
\end{lemma}

\paragraph{Proof.}
Suppose first that the parameter set contains a factor with $s=p^n > \sqrt{N}$ levels.
If $m$ is even then a geometric $OA(p^m , (p^{m-n} )^{p^n} (p^n)^1)$ exists
by Section 5, and $p^{m-n}$ is the largest number of levels other factors can have if there is an $s$-level factor.
Since there are $p^n$ factors with $p^{m-n}$ levels, the existence of any array with one
$s$-level factor and at most $p^{m/2}$ factors with $\le p^{m-n}$ levels follows immediately.
The case that $m$ is odd follows similarly.

We now assume that all $s_i \le \sqrt{N}$.
If $m$ is even, $N=p^{2r}$, then a geometric
\beql{Eq61}
OA (p^{2r} , (p^r )^{ p^r +1})
\eeq
exists by Section 5.
Any parameter set with all $s_i \le p^r$ and $k \le p^r +1$ is dominated by \eqn{Eq61} and so is realized.
If $m$ is odd, $N= p^{2r+1}$, then a geometric
\beql{Eq62}
OA(p^{2r+1}, (p^r)^{p^{r+1} +1} )
\eeq
also exists by Section 5.
Any parameter set with all $s_i \le p^r$ and $k \le p^{r+1} +1$ is dominated by \eqn{Eq62} and so is also realized.~~~$\bsq$

Since the number of factors in a parameter set is less than or equal to the number of degrees of freedom
\eqn{Eq61a}, Lemma \ref{le6} immediately implies that any
parameter set
$(N=p^m, s_1^{k_1} s_2^{k_2} \ldots )$ with at most
$p^{\lfloor (m+1)/2 \rfloor} +1$ degrees of freedom is realized by an orthogonal array.
However, Theorem \ref{th4} is much stronger.

\paragraph{Proof of Theorem \ref{th4}.}
We will show that any parameter set
$(p^m, p^{k_1} (p^2)^{k_2} (p^3)^{k_3} \ldots )$ satisfying
\beql{Eq40}
\sum_{i \ge 1} k_i (p^i -1) \le p^{3m/4}
\eeq
is realized by a geometric orthogonal array, where $p$ is any prime.
To simplify the notation we assume $m=4r$
is a multiple of 4.
The arguments in the other three cases require only minor modifications and are left to the reader.

From \eqn{Eq40} we have
\beql{Eq41}
k_{r+1} + k_{r+2} + \cdots + k_{4r-1} \le p^{2r} +1 ~,
\eeq
and so by Lemma \ref{le6} a geometric
$$OA( p^{4r} , (p^{r+1})^{k_{r+1}} (p^{r+2})^{k_{r+2}} \ldots
(p^{4r-1} )^{k_{4r-1}} )
$$
exists.

We now proceed by induction.
Let $H_n$ be the hypothesis that every parameter set
\beql{Eq50}
(p^{4r} , p^{b_1} (p^2)^{b_2} \ldots
(p^r)^{b_r} (p^{r+1} )^{k_{r+1}} \ldots
(p^{4r-1} )^{k_{4r-1}} )
\eeq
with
$b_i \le k_i$ for $1 \le i \le r$ and
$$b_1 + b_2 + \cdots + b_r = n$$
can be realized by a geometric orthogonal array constructed using disjoint
subspaces of $PG(4r-1, p)$.
We have shown that $H_0$ holds.
Suppose $H_n$ holds with
$$n < k_1 + k_2 + \cdots + k_r ~.$$
We will show that we can increase $b_r$ by 1 and still realize the
parameter set, thus establishing $H_{n+1}$.

To show this, note that the number of projectively
distinct nonzero points in all the subspaces in \eqn{Eq50} is at most
$$\sum_{i=1}^{4r-1} k_i \frac{p^i -1}{p-1} -1 \le \frac{p^{3r}}{p-1} -1 ~.$$
However, by Theorem 1 of (Thas, 1995), Section 7,
if a subset of $PG (4r-1, p)$ contains fewer than
$$(p^{3r+1} -1)/ (p-1) ~~\mbox{points} ~,
$$
there is a subspace $PG(r-1, p)$ disjoint from it.
Since
$$\frac{p^{3r}}{p-1} -1 < \frac{p^{3r+1}-1}{p-1} ~,$$
such a subspace exists and we can use it to augment $b_r$ by 1.

By induction, we can realize the parameter set $(p^{4r}, p^{k_1} (p^2)^{k_2} \ldots (p^{4r-1})^{k_{4r-1}} )$, as required.~~~$\bsq$

\paragraph{Proof of lower bound of Theorem \ref{th2}.}
First suppose $N = 2^{2r}$.
From Theorem \ref{th4}, every parameter set
$$(2^{2r} , 2^{k_1} 4^{k_2} 8^{k_3} \ldots (2^r )^{k_r} )$$
with at most $2^{3r/2}$ degrees of freedom can be realized.
The lower bound of Theorem \ref{th2} now follows in the same way that
we proved the upper bound in Section 4.
If $N = 2^{2r+1}$ we use the previous case
together with $T(N) \ge T(N/2)$.~~~$\bsq$

It would be nice to have analogues of Theorem \ref{th4} and Lemma \ref{le6} when $N$ is not a prime power!

\subsection*{Acknowledgments}
We thank Michele Colgan for computing the properties of the lattices
$\La'_N$ shown in Table \ref{ta3}.
The research of John Stufken
was supported by NSF grant DMS-9803684.

%

\clearpage
\section*{References}

\begin{description}

\item[~~~]
Bosma, W. and Cannon, J. (1995).
{\em Handbook of Magma Functions},
Sydney.
 
\item[~~~]
Bosma, W., Cannon, J. and Mathews, G. (1994).
Programming with algebraic structures:
Design of the Magma language, in
{\em Proceedings of the 1994 International Symposium on Symbolic and Algebraic Computation},
M. Giesbrecht, Ed.,
Association for Computing Machinery, 52--57.
 
\item[~~~]
Bosma, W., Cannon, J.  and Playoust, C. (1997).
The Magma algebra system I:
The user language.
{\em J. Symb. Comp.},
{\bf 24}, 235--265.

\item[~~~]
Calderbank, A. R., Rains, E. M., Shor, P. W. and Sloane, N. J. A. (1998).
\htmladdnormallink{Quantum error correction via codes over $GF(4)$}{http://www.research.att.com/~njas/doc/qc2.pdf}.
{\em IEEE Trans. Information Theory}, {\bf 44}, 1369--1387.

\item[~~~]
Eisfeld, J., Storme, L. and Sziklai, P. (1999).
Maximal partial spreads in finite projective spaces, preprint.

\item[~~~]
Hedayat, A. S., Sloane, N. J. A. and Stufken, J. (1999).
{\em \htmladdnormallink{Orthogonal Arrays: Theory and Applications}{http://www.research.att.com/~njas/doc/OA.html}.}
New York: Springer-Verlag.

\item[~~~]
Hoffman, K. and Kunze, R. (1961).
{\em Linear Algebra},
Engelewood Cliffs, NJ: Prentice-Hall.

\item[~~~]
Kimura, H. (1994a).
Classification of Hadamard matrices of order 28 with Hall sets.
{\em Discrete Math.},
{\bf 128},
257--268.

\item[~~~]
Kimura, H. (1994b).
Classification of Hadamard matrices of order 28.
{\em Discrete Math.},
{\bf 133},
171--180.

\item[~~~]
Sloane, N. J. A. (1999).
The On-Line Encyclopedia of Integer Sequences.
Published electronically at \htmladdnormallink{http://www.research.att.com/$\sim$njas/sequences/}{http://www.research.att.com/~njas/sequences/}.

\item[~~~]
Sloane, N. J. A. and Stufken, J. (1996).
\htmladdnormallink{A linear programming bound for orthogonal arrays with mixed levels}{http://www.research.att.com/~njas/doc/mixed.pdf}.
{\em J. Statist. Plann. Infer.},
{\bf 56}, 295--305.

\item[~~~]
Thas, J. A. (1995).
Projective geometries over a finite field,
Chapter 7 of {\em Handbook of Incidence Geometry}, edited by
Buekenhout, F. Amsterdam: North-Holland.

\item[~~~]
Trotter, W. F. (1995).
Partially ordered sets,
Chapter 8 of {\em Handbook of Combinatorics}, edited by Graham, R. L.,
Gr\"{o}tschel, M. and Lov\'{a}sz, L. Amsterdam: North-Holland;
Cambridge, MA: MIT Press.

\item[~~~]
Welsh, D. J. A. (1976).
{\em Matroid Theory.}
London: Academic Press.

\item[~~~]
Wu, C. F. J. (1989).
Construction of $2^m 4^n$ designs via a grouping scheme.
{\em Annals Statist.}, {\bf 17}, 1880--1885.
\end{description}

\end{document}